\documentclass[11pt]{article}

\usepackage{amsmath,amsfonts,amscd, amssymb,theorem}

\newcommand{\doubleheaddownarrow}{\big\downarrow\kern-3.325mm\downarrow}
\newcommand\Oh{\mathcal O}
\newcommand\Z{\mathbb Z}

\newcommand\om{\omega}

\newcommand\Hom{\operatorname{Hom}}

\newcommand{\skewentry}{ \kern -1cm -\mathrm{sym} \kern -1cm}
\newcommand\Spec{\operatorname{Spec}}

 \newtheorem{theorem}{Theorem}[section]
 \newtheorem{lemma}[theorem]{Lemma}
 \newtheorem{prop}[theorem]{Proposition}
 
 {
 \theorembodyfont{\rmfamily}
 \newtheorem{defn}[theorem]{Definition}

 \newtheorem{rem}[theorem]{Remark}

 }
 \newenvironment{pf}{\paragraph{Proof}}{\par\medskip}
 
 \newcommand{\qed}{\ifhmode\unskip\nobreak\fi\quad\ensuremath\square}
\newcommand{\QED}{\ifhmode\unskip\nobreak\fi\quad\ensuremath{\mathrm{QED}}}

\numberwithin{equation}{section}

\title{ Remarks on Type III Unprojection }

\author{Stavros Argyrios Papadakis\thanks{email: spapad@maths.warwick.ac.uk}}
\date {January 2005}

\begin{document}
\maketitle

\begin {abstract}
Type III unprojection plays a very important role in the birational
geometry of Fano threefolds (cf. \cite{CPR}, \cite{Ki}, \cite{BZ}). 
According to \cite{Ki} p.~43, it was first introduced by A. Corti
on his calculations of Fano threefolds of genus 6 and 7.

It seems that at present a general definition of type III unprojection
is still missing. After proving in Section~\ref{sec!residual}
some general facts about residual ideals, we propose a definition
for the generic Type III unprojection (Definition~\ref{dfn!cimaindfn}),
and prove in Theorem~\ref{thm!basicthrm} that it gives a Gorenstein
ring. 
\end{abstract}

\section {Introduction}  \label{sec!intro}

Unprojection is a philosophy, which aims to construct
and analyse Gorenstein rings in terms of simpler ones. Geometrically,
it can be considered as an inverse of projection,
and as a method for constructing birational 
contractions.

So far they have appeared at least four types of unprojection.
Unprojection of type Kustin--Miller (or type I) 
(\cite{KM},\cite{PR},\cite{P}), type II (\cite{CPR}, \cite{P2}),
type IV (\cite{R}), and type III which is the subject of
the present work.

Type III unprojection is the residual unprojection of the
simplest unprojection of type Kustin--Miller, the one of
a complete intersection inside a complete intersection 
(\cite{P} Section~4). In geometry, it appears 
as the contraction to the 'other direction' from the
middle point of a Sarkisov link between Fano threefolds. 
A geometric example is treated in \cite{Ki} Examples~4.6 and 9.16.

The structure of the paper is as follows. In Section~\ref{sec!residual} 
we define the residual ideal of an unprojection of type 
Kustin--Miller (Definition~\ref{dfn!residdfn1}) and 
prove some basic properties of it. In 
Definition~\ref{dfn!cimaindfn}
we propose a definition for the generic type III
unprojection.  An important
aim in the theory of unprojections is to get a Gorenstein
unprojection ring, and in Theorem~\ref{thm!basicthrm},
which is our main result,  we
give a presentation of the unprojection ring, from which
the Gorensteiness follows immediately. As mentioned
in Remark~\ref{rem!gencase} it will be nice to have
a treatment of the general Type III unprojection.

\cite{Ki} contains more discussion about the various unprojections and 
their applications to algebraic geometry.

\section {Residual ideal of unprojection} \label{sec!residual}

Assume $I_D \subset \Oh_X$ is a pair of an ideal and a ring satisfying the
assumptions of \cite{PR} Section~1. That is,  $\Oh_X$ is a 
local commutative Gorenstein ring, and $I_D$ is a codimension one ideal 
with  $\Oh_X/I_D$ Gorenstein. Moreover, we fix   an 
$\Oh_X$-regular element $z \in I_D$, such that $(z) \not= I_D$.

\begin {defn}   \label {dfn!residdfn1}
The residual ideal $I_r$ of $I_D$ with respect to the element 
$z$ is the ideal
\[
   I_r = \{ f(z) \colon f \in \Hom_{\Oh_X}(I_D, \Oh_X) \} \subset \Oh_X.
\]
\end{defn}

\paragraph {} We set $\Oh_D = \Oh_X/I_D$, $\Oh_r = \Oh_X/I_r$.

\begin{rem}   \label {rem!twogens}
By the adjunction sequence of \cite{PR} Section~1,  there exists
$f \in \Hom_{\Oh_X}(I_D, \Oh_X)$ such that 
\[
    I_r = (z, f(z)).
\]
\end{rem}

\begin {lemma}   \label{lem!codimisone}
The codimension of $I_r$ in $\Oh_X$ is one.
\end{lemma}

\begin{pf}
It is enough to notice that $z, f(z)$ is not a regular
$\Oh_X$-sequence, since for $w \in I_D \setminus (z)$
we have
\[
        w f(z) - z f(w) =0.
\] 
\QED \medskip
\end{pf}

\paragraph{} From the definition of $I_r$ we have 
\begin {equation}   \label{eq!colonforIr}
     I_r = (z) \colon I_D. 
\end{equation}
Indeed, let  $f(z) \in I_r$ for some 
$f \in \Hom_{\Oh_X}(I_D, \Oh_X)$. Then, for any
$a \in I_D$ we have
\[ 
      a f(z) = f(az) = z f(a) \in (z).
\]
Conversely, if $w \in (z) \colon I_D$, multiplication
by $w/z$ is an element of  $\Hom_{\Oh_X}(I_D, \Oh_X)$
and $(w/z)z = w$, so $w \in I_r$.

Since $\Oh_D$ is Gorenstein we have more.

\begin{prop}   \label{prop!residprop} 
The ring $\Oh_r$ is Cohen--Macaulay, and 
\[  
         I_D = (z) \colon I_r.
\]
\end {prop}

\begin{pf}
We work over the ring $\Oh_X/(z)$. Using 
(\ref{eq!colonforIr}) and  \cite{Ei}~Theorem~21.23, 
we  get that $\Oh_r$ is Cohen--Macaulay
and 
\[
    I_D/(z) = 0 \colon I_r/(z), 
\]
which implies $ I_D = (z) \colon I_r$.
\QED \medskip
\end{pf}           

Denote by $\om_r$ the dualizing module of the local ring $\Oh_r$, by
\[
   \phi_1  \colon \Oh_X \to (z) 
\]
the map with  $ \phi_1(a)=az$, and by 
\[
   \phi_2  \colon   \Hom_{\Oh_X}(I_r,\Oh_X) \to I_D
\]
the map with
\[
   \phi_2 (f) = f(z).
\]

We have a commutative diagram 
\[
 \renewcommand{\arraystretch}{1.5}
 \begin{matrix}
 0 &\to&  \Oh_X   &\to&   \Hom_{\Oh_X}(I_r,\Oh_X) 
                  &\to&\om_r&\to&0                  \\
 &&  \downarrow{\phi_1} && \downarrow{\phi_2} &&    
                  \downarrow{\phi_3}     \\
 0 &\to& (z) &\to& I_D &\to&  I_D/(z) &\to&0 
 \end{matrix}
 \]
where the first row is the adjunction sequence of
\cite{PR} Section~1,
and $ \phi_3$ is the induced map. 

The following proposition follows immediately from 
Proposition~\ref{prop!residprop}.
 
\begin{prop}    \label{prop!threeiso}
 The vertical maps $\phi_i$, $i=1,\dots ,3$ in the above 
diagram are isomorphisms of $\Oh_X$-modules.
\end{prop}

\section {Generic type III unprojection}  \label{sec!generic}

Fix $n \geq 1$. Let $\Oh_{amb}$ be the polynomial ring  
$\Oh_{amb} = \Z [a_{ij}, z_j]$, with $1 \leq i \leq n, \ 
1 \leq j \leq n+1$, $M$ the $n \times n+1$ matrix
\[
   M =  \begin {pmatrix}
             a_{11} & \dots & a_{1,n+1}   \\
             \vdots &  \vdots   &  \vdots \\
             a_{n1} & \dots & a_{n,n+1} 
     \end {pmatrix}, 
\] 
and  $\Delta_1$  the determinant of the submatrix
of $M$ obtained by deleting the first column. We also set   
\[
 f_i(z) = \sum_{j=1}^{n+1} a_{ij}z_j,
\]
for $ 1 \leq i \leq n $. Denote by $I_X$ the ideal
\[
   I_X = (f_1, \dots ,f_n) \subset \Oh_{amb},
\]
by $\Oh_X = \Oh_{amb}/I_X $ the quotient ring, and by $I_D$
the codimension ideal $I_D =(z_1, \dots ,z_{n+1})$
of $\Oh_X$. Since both $z_1, \dots ,z_{n+1}$ and $f_1, \dots
f_n$ are regular sequences of $\Oh_{amb}$, we have 
that the ideal $I_D$ of $\Oh_X$  has codimension one. 

\begin{prop} \label{prop!number007}
The residual ideal  (in the sense of Definition~\ref{dfn!residdfn1})
of $I_D \subset \Oh_X$ with respect to the element $z_1$ is the ideal 
\[
          I_r = ( z_1, \Delta_1) \subset \Oh_X.
\]
\end{prop}

\begin {pf}
 It follows immediately from \cite{P}~Theorem~4.3.
\QED \medskip
\end{pf}

\begin{prop}  \label {prop!inters}
   We have 
\begin {equation} \label{eqn!inters}
    (z_1) = (z_1, \dots ,z_{n+1}) \cap (z_1,\Delta_1) \subset \Oh_X
\end{equation}
 (equality of ideals of $\Oh_X$).  
Moreover, the quotient ring  $\Oh_X/ (z_1, \Delta_1)$
 is a Cohen--Macaulay  integral domain.
\end{prop}

\begin{pf}

To prove (\ref{eqn!inters}) it is enough to prove the equality 
\[
    (z_1, f_1, \dots , f_n) = (z_1, \dots ,z_{n+1}) \cap
                               (z_1, f_1, \dots ,f_n, \Delta_1),
\]
of ideals of $\Oh_{amb}$.

     By Cramer's rule (cf. \cite{P} Lemma~4.2)
 we have that $ \Delta_1 z_j \in (z_1, f_1, 
\dots, f_n)$ for $1 \leq j \leq n+1$. Let 
\[
   h = c \Delta_1 \in (z_1, \dots ,z_{n+1}) \cap (\Delta_1),
\]
for an element $c \in \Oh_{amb}$. Then $c \in  (z_1, \dots ,z_{n+1})$,
therefore
\[
    (z_1, f_1, \dots , f_n) = (z_1, \dots ,z_{n+1}) \cap
                               (z_1, f_1, \dots ,f_n, \Delta_1).
\]

The fact that the quotient ring  $\Oh_{amb}/ (z_1, f_1, \dots ,
                                              f_n, \Delta_1)$
is a Cohen--Macaulay integral domain  follows 
as in \cite{H} Section~3 Example~4. The  Cohen-Macaulayness also
follows from Proposition~\ref{prop!residprop}. 
\QED \medskip
\end{pf}

We denote by $K(X)$ the quotient field of the integral domain $\Oh_X$, 
and we set 
\[
     I_r^{-1} = \{ f \in K(X) \colon f I_r \subset \Oh_X \}.
\]

\begin{defn}   \label {dfn!cimaindfn}
The generic type III unprojection of $I_D \subset \Oh_X$ with respect 
to the element $z_1$ is the $\Oh_X$-subalgebra of $K(X)$
\[
         \Oh_X [I_r^{-1}] \subset K(X)
\]
generated by $I_r^{-1}$.
\end {defn}

\begin {rem}  We are extending $\Oh_X$ by including rational
functions with denominators in $I_r$. Geometrically, it corresponds 
to 'contracting' the codimension one subscheme of 
$X = \Spec \Oh_X$ defined by the ideal $I_r$. It is 
interesting to compare with the proof of Castelnuovo's 
contractibility criterion given in \cite{Be} p.~20.
\end{rem}

\paragraph{} Using Proposition~\ref{prop!threeiso}, we see that 
$\Oh_X [I_r^{-1}]$ is generated as an $\Oh_X$-algebra by 
$s_2, \dots ,s_{n+1} \in K(X)$, with
\[
    s_i = \frac{z_i}{z_1}. 
\]

Define the ring homomorphism 
\[
     \phi \colon \Oh_{amb}[T_2, \dots ,T_{n+1}] \to K(X),
\]
which restricted to $\Oh_{amb}$ is the natural projection 
to $\Oh_X$, and
\[
     \phi (T_i) = s_i,
\]
for $1 \leq i \leq n+1$. Also define 
 polynomials $f_i(T) \in  \Oh_{amb}[T_1, \dots ,T_{n+1}]$ 
by
\[   
   f_i(T_1, \dots ,T_n) = \sum_{j=1}^{n+1} a_{ij}T_j,
\] 
for $1 \leq i \leq n$.

Our aim is to prove the following theorem.

\begin{theorem}  \label{thm!basicthrm}
We have
\begin {equation} \label {eqn!basiceqns}
  \ker \phi = (z_i -z_1T_i, f_j(1,T_2, \dots ,T_{n+1})),
\end {equation}
with indices $2 \leq i \leq n+1$ and $1 \leq j \leq n$. 
As a consequence,   $\Oh_X [I_r^{-1}]$ is isomorphic
to a polynomial ring, hence it is Gorenstein.
\end {theorem}

For this, we will use the ring homomorhism 
\[
     \psi \colon \Oh_{amb}[T_1, T_2, \dots ,T_{n+1}] \to \Oh_X,
\]
which restricted to $\Oh_{amb}$ is the natural projection 
to $\Oh_X$, and
\[
     \psi (T_i) = z_i.
\]
for $1 \leq i \leq n+1$.

\begin {rem} We remark that the equations obtained 
in \cite{Ki} Example~9.16 are a specialization of 
those in (\ref{eqn!basiceqns}).
\end {rem}

\paragraph{}The following proposition is immediate.

\begin{prop} \label{prop!kernelofpsi}
We have
\[
   \ker \psi = (T_i-z_i, f_j(z)),
\]
with indices $1 \leq i \leq n+1$ and $1 \leq j \leq n+1$. 
\end {prop}

\paragraph {}
On the polynomial ring $\Oh_{amb}[T_1, \dots ,T_{n+1}]$, we put 
weight $0$ for elements of $\Oh_{amb}$ and $1$ for the $T_i$.

\begin{prop} \label{prop!prophomog}
Denote by 
\[
    J \subset \ker \psi
\]
the biggest homogeneous ideal of   $\Oh_{amb}[T_1, \dots ,T_{n+1}]$
contained in $\ker \psi$.  Then
\[
    J = (z_iT_j-z_jT_i, f_p(z), f_p(T)),
\]
with indices $1 \leq i,j \leq n+1$ and $1 \leq p \leq n$.
\end{prop}

\begin {pf}

\paragraph {}  During the proof of the  Proposition~\ref{prop!prophomog}
we use the following notations.

We set  
\[
    R_1 = \Oh_{amb}[T_1, T_2, \dots ,T_{n+1}],
\]
also
\[ 
  z = (z_1, \dots ,z_{n+1}),\ T = (T_1, \dots ,T_{n+1}),\ \
\]
and
\[
   f(z) = (f_1(z), \dots, f_n(z)), \  \
   f(T) = (f_1(T), \dots, f_n(T)).
\]

For any $i \geq 0$, by
\[
   c_i = (c_i^1, \dots , c_i^{n+1})
\]
we will denote an $n+1$-tuple of elements $c_i^j \in R_1$ 
and similarly, by 
\[
  l_i = (l_i^1, \dots , l_i^n)
\]
we will denote an $n$-tuple of elements 
$l_i^j \in R_1$.

For two elements $a=(a_1, \dots ,a_p), 
b=(b_1, \dots ,b_p)$ of $(R_1)^p$, 
we set
\[
    <a,b> = \sum_{i=1}^{p}a_ib_i.
\]

\paragraph {Claim 1} 
We have 
\[
     z_k f_l(T) -  T_k f_l(z) \in (z_i T_j - z_j T_i).
\] 

Indeed 
\[
  z_k f_l(T) -  T_k f_l(z)= \sum_j  a_{lj}(z_kT_j-T_kz_j).
\]

\paragraph {Claim 2} Let
\[
   h_0 = <-c_0, z> + <l_0, f(z)>,
\]
and for $1 \leq i \leq s$,
\[
   h_i = <-c_i, z> + <c_{i-1},T>+  <l_i, f(z)>,
\]
for elements $c_i \in (R_1)^{n+1}, \ l_i \in (R_1)^n$. 

Assume $h_0= \dots =h_s =0$. Then
\[
   <c_s, T> \ \in  J.
\]

\paragraph{}  Proof of Claim 2 by induction on $s$. 

Assume $s=0$. From $h_0 = 0$ it follows that
\[
   \sum_{i=1}^{n+1}  z_i m_i  = 0.
\]
with 
\[
   m_i = c_0^i+ \sum_{j=1}^n l_0^ja_{ji},
\]
for $1 \leq i \leq n+1$.

Since the $z_i$ form a regular sequence, we have 
that   $(m_1, \dots ,m_{n+1})$ is a linear combination
of Koszul relations $K_{ij}=(b_1, \dots ,b_{n+1})$, where 
$b_t = 0$ for $t \notin \{i,j\}$, $b_i = z_j$ and $b_j = -z_i$. 
As a consequence, $<c_0, T> \ \in J$.

Assume now $h_0 = h_1 = \dots = h_s = 0$ for $s \geq 1$. 
By the case $s=0$  we have
\[
   h_1 = <-c_1,z> + u_1 + u_2 + u_3,
\]
where $u_1$ is a linear combination of terms $f_i(z)$, $u_2$ is a 
linear combination of terms $f_i(T)$, and $u_3$ is a linear
combination of term $z_iT_j-z_jT_i$. Since no $z_i$ appears in 
$f_i(T)$,  from $h_1 =0$ it follows that we can assume
that each coefficient of $f_i(T)$ is in $(z_1,\dots ,z_{n+1})$. Using 
Claim 1, we can assume $u_2 = 0$ by adding terms to $u_1$ and $u_3$.

Now, since $u_3$ is a linear combination  of terms $z_iT_j-z_jT_i$,
we can change $c_1$ to $c'_1$ in such a way that 
 $h_1 = <c'_1,z> + \sum  \text{coeff} f_i(z)$
and $<c'_1,T> = <c_1, T>$. Claim 2 follows from the 
inductive hypothesis, since $h_1$ has the same form as $h_0$.

\paragraph {} The same arguments prove the following.

\paragraph {Claim 3}  Let
\[
   h_0 = <c_0, T> + <l_0,f(z)>,
\]
and for $1 \leq i \leq s$
\[
   h_i = <c_i, T> + <-c_{i-1},z>+  <l_i, f(z)>.
\]
Assume $h_0= \dots =h_s =0$. Then
\[
   <c_s, z> \ \in  J.
\]

\paragraph {} We now finish the proof of 
Proposition~\ref{prop!prophomog}.

Let $h \in  (T_i-z_i, f_j(z))$ be a homogenous
element.
Write
\[
   h = \sum_{p=1}^{n+1} c^p (T_p-z_p) + \sum_{p=1}^n l^pf_p(z),
\]
with 
\[
     c^p = c_0^p+ \dots +c_s^p, \quad 
     l^p = l_0^p+ \dots +l_{s+1}^p,
\]
and each $c^p_i, l^p_i$ homogeneous of degree $i$.
Then
\[
     h = h_0 + \dots + h_{s+1}
\]
is the decomposition of $h$ into homogeneous components, with
\[
         h_0 = <-c_0,z> + <l_0,f(z)>,
\]
and for $1 \leq i \leq s$
\[
      h_i = <-c_i, z> + <c_{i-1},T> + <l_i, f(z)>
\]
and 
\[
    h_{s+1} =  < c_s, T> + <l_{s+1}, f(z)>.
\]
Since $h$ is assumed homogeneous, exactly one $h_q \not= 0$. The result
follows from Claims 2 and 3.
\QED \medskip
\end{pf}

We now give the proof of Theorem~\ref{thm!basicthrm}.

\paragraph  {Proof of Theorem~\ref{thm!basicthrm} \newline} 

Assume  $h(T_2,\dots ,T_{n+1}) \in \ker \phi$,
with total degree (with respect to $T_i$) equal to $d$. 
Set
\[
    h'(T_1, \dots ,T_{n+1}) = T_1^dh(\frac{T_2}{T_1},
    \dots , \frac{T_{n+1}}{T_1}) \in \Oh_{amb} [T_1, \dots ,T_n].
\]
We have $h' \in \ker \psi$ homogeneous, and the Theorem follows 
from  Proposition~\ref{prop!prophomog}.
\QED \medskip

\begin {rem} \label{rem!gencase}
A very interesting open question is to generalise 
Definition~\ref{dfn!cimaindfn} and Theorem~\ref{thm!basicthrm} so
as to cover the general case of Type III unprojection.
\end{rem}

\paragraph {ACNOWLEDGMENTS}  I wish to thank 
Miles Reid, Frank--Olaf Schreyer and Francesco Zucconi for many
useful conversations. This work was financially supported by
the Deutschen Forschungsgemeinschaft Schr 307/4-2.

\bigskip
\noindent
Stavros Papadakis, \\
Mathematik und Informatik,  Geb. 27 \\
Universitaet des Saarlandes \\
D-66123 Saarbruecken, Germany \\
e-mail: spapad@maths.warwick.ac.uk


\begin {thebibliography} {xxx}



\bibitem[Be]{Be} 
Beauville A.,
\textsl{Complex algebraic surfaces, Second Edition}.
London Mathematical Society Student Texts 34,
CUP 1996

\bibitem[BZ]{BZ}
Brown G. and Zucconi F., work in progress

\bibitem[CPR]{CPR}
Corti A., Pukhlikov A. and Reid M.,
\textsl{Birationally
rigid Fano hypersurfaces}, in Explicit birational geometry
of 3-folds,
A. Corti and M. Reid (eds.), CUP 2000, 175--258

\bibitem[Ei]{Ei}
Eisenbud D., \textsl {
Commutative algebra,
with a view toward algebraic geometry}.
Graduate Texts in Mathematics 150,
Springer--Verlag 1995

\bibitem[H]{H}
Herzog J.,
\textsl{
Certain complexes associated to a sequence and a matrix},
Manuscripta Math. {\bf  12} (1974) 217--248


\bibitem[KM]{KM} Kustin A. and Miller M., \textsl {
Constructing big Gorenstein ideals from small ones}.
J. Algebra {\bf 85} (1983),  303--322

\bibitem[P]{P} Papadakis S., \textsl {
Kustin--Miller unprojection with complexes},
J. Algebraic Geometry {\bf 13}  (2004), 249-268

\bibitem[P2]{P2} Papadakis S.,\textsl {
Type II unprojection}, submitted, 17~pp.

\bibitem[PR]{PR} Papadakis  S. and Reid M., \textsl {
Kustin--Miller unprojection without complexes},
J. Algebraic Geometry {\bf 13}  (2004), 563-577

\bibitem[R]{R}
Reid M.,
\textsl{Examples of type IV unprojection},
math.AG/0108037, 16~pp.

 \bibitem[Ki]{Ki} Reid M., \textsl {
Graded Rings and Birational Geometry},
in Proc. of algebraic symposium (Kinosaki, Oct 2000),
K. Ohno (Ed.) 1--72, available from
www.maths.warwick.ac.uk/$\sim$miles/3folds



\end{thebibliography}
\end{document}